\documentclass[12pt,twoside]{article}
\title{The primitive derivation and freeness of multi-Coxeter arrangements}
\author{Masahiko Yoshinaga \thanks{RIMS, Kyoto university Kyoto, 
606-8502 Japan, 
e-mail address: yosinaga@kurims.kyoto-u.ac.jp}}
\date{May 5, 2002}

\usepackage{amssymb}

\begin{document}

\maketitle

\newtheorem{Def}{Definition}
\newtheorem{Prop}[Def]{Propositon}
\newtheorem{Thm}[Def]{Theorem}
\newtheorem{Lemma}[Def]{Lemma}
\newtheorem{Cor}[Def]{Corollary}
\newtheorem{Rem}[Def]{Remark}
\newtheorem{Problem}[Def]{Problem}

\begin{abstract}
We will prove the freeness of multi-Coxeter arrangements 
by constructing a basis of the module of vector fields which 
contact to each reflecting hyperplanes with some multiplicities 
using K. Saito's theory of primitive derivation. 



\end{abstract}

\section{Introduction}\label{sec:intro}

Let $V$ be a Euclidean space over $\mathbb{R}$ with finite dimension $\ell$ 
and inner product $I$. Let  
$W\subset \mathrm{O}(V,I)$ be a finite irreducible reflection group and 
$\mathcal{A}$ the corresponding Coxeter arrangement i.e. the 
collection of all reflecting hyperplanes of $W$. 
For each $H\in \mathcal{A}$, we fix a defining equation 
$\alpha_H \in V^\ast$ of $H$. 

In \cite{ter-multi}, H. Terao constructed a free basis of 
$\mathbb{R}[V]$-module 
\begin{equation}
\mathrm{D}^m (\mathcal{A}):=\{\delta \in \mathrm{Der}_V\ |\ 
\delta \alpha_H \in (\alpha_H^m),\ \forall H\in \mathcal{A}\ \}, 
\end{equation}
($m\in \mathbb{Z}_{\geq 0}$). 
The purpose of this paper is to 
construct a basis by a simpler way 
using Saito's result and give a generalization. 

For given multiplicity 
$\tilde{m}:\mathcal{A}\rightarrow \mathbb{Z}_{\geq 0}$, 
we say that the multi-Coxeter 
arrangement $\mathcal{A}^{(\tilde{m})}$ is free if the module 

\begin{equation}
{\rm D}(\mathcal{A}^{(\tilde{m})}):=\{ \delta \in \mathrm{Der}_V\ |\ 
\delta \alpha_H \in (\alpha_H^{\tilde{m}(H)}),\ \forall H\in \mathcal{A}\ \}
\end{equation}
is a free $\mathbb{R}[V]$-module \cite{zie}. 
Then our main result is 
\begin{Thm}\label{thm:main}
Let $\tilde{m}$ be a multiplicity satisfying 
$\tilde{m}(H)\in \{0,1\}$ for all $H\in \mathcal{A}$. Suppose 
the multi-Coxeter 
arrangement $\mathcal{A}^{(\tilde{m})}$ is free, then 
$\mathcal{A}^{(\tilde{m}+2k)}\ (k\in \mathbb{Z}_{\geq 0})$ is also free, 
where the new multiplicity $\tilde{m}+2k$ take value $\tilde{m}(H)+2k$ at 
$H\in \mathcal{A}$.
\hfill$\square$
\end{Thm}
We construct a basis in Theorem \ref{thm:bas}. 

\noindent
We note that $\mathcal{A}^{(\tilde{m})}$ is not necesarily free for 
$\tilde{m}:\mathcal{A}\rightarrow \{0,1\}$. 
If we apply Theorem \ref{thm:main} for 
$\tilde{m}(H)\equiv 0$ or $\tilde{m}(H)\equiv 1$, 
we obtain the freeness of $\mathrm{D}^{2k} (\mathcal{A})$ or
$\mathrm{D}^{2k+1} (\mathcal{A})$. Terao's basis 
is expected to coincide with that of ours. 

The original motivation to study the module 
$\mathrm{D}(\mathcal{A}^{(\tilde{m})})$ came from the study of 
structures of the relative de Rham cohomology 
$\mathrm{H}^\ast(\Omega_{\mathfrak{g}/S}^\bullet)$ 
of the adjoint quotient map 
$\chi:\mathfrak{g}\rightarrow S:=\mathfrak{g}//\mathrm{ad}(G)$
of a simple Lie algebra $\mathfrak{g}$. 
In the case of $ADE$ type Lie algebras, an isomorphism as 
$\mathbb{C}[S](=\mathbb{C}[\mathfrak{g}]^G=\mathbb{C}[\mathfrak{h}]^W)$-modules
(where $\mathfrak{h}$ is a Cartan subalgebra) 
\begin{displaymath}
\mathrm{H}^2(\Omega_{\mathfrak{g}/S}^\bullet)
\cong
\mathrm{D}^5(\mathcal{A})^W
\end{displaymath}
is obtained \cite{yos}. 

But for $BCFG$ type Lie algebras, 
because the $W$ action on $\mathcal{A}$ is not transitive, 
$\mathrm{H}^2(\Omega_{\mathfrak{g}/S}^\bullet)$ is expected to be 
isomorphic to 
the module $\mathrm{D}(\mathcal{A}^{(\tilde{m})})^W$ with a suitable 
multiplicity $\tilde{m}:\mathcal{A}\rightarrow \mathbb{Z}_{\geq 0}$ 
which is not constant.

\section{K. Saito's results on primitive derivation}\label{sec:saito}

In this section, we fix  notations and recall some results. 

Let $x_1,\cdots ,x_\ell \in V^\ast$ be a basis of $V^\ast$ and 
$P_1, P_2,\cdots ,P_\ell \in \mathbb{R}[V]^W$ be the homogeneous 
generators of $W$-invariant polynomials on $V$ such that 
$\mathbb{R}[V]^W=\mathbb{R}[P_1, P_2,\cdots ,P_\ell]$ with 
\begin{displaymath}
\deg P_1 \leq \deg P_2 \leq \cdots \leq \deg P_\ell=:h. 
\end{displaymath}
Then it is classically known \cite{bou} that 
\begin{equation}\label{eq:num}
|\mathcal{A}|=\frac{h \ell }{2}
\end{equation}
and 
\begin{equation}\label{eq:deg}
\deg P_{\ell-1} < h.
\end{equation}
It follows from (\ref{eq:deg}) that 
the rational vector field (with pole along $\bigcup_{H\in \mathcal{A}}H$)
$D:=\frac{\partial}{\partial P_\ell}$ 
on $V$ is uniquely determined up to non-zero constant factor 
independently on 
the generators $P_1,\cdots,P_\ell$. 
We call $D$ the primitive vector field. If we fix generators 
$P_1,\cdots,P_\ell$, then 
$\frac{\partial}{\partial P_1},\cdots,\frac{\partial}{\partial P_{\ell -1}}$ 
are able to be considered as rational vector fields on $V$. 
Since the Jacobian is 
\begin{displaymath}
Q:=\prod_{H\in \mathcal{A}}\alpha_H
\dot{=}
\frac{\partial (P_1,\cdots ,P_\ell)}{\partial (x_1,\cdots ,x_\ell)}, 
\end{displaymath}
$D$ is symbolically expressed as
\begin{displaymath}
D\dot{=}\frac{1}{Q}
\det
\left(
\begin{array}{cccc}
\frac{\partial P_1}{\partial x_1} & \cdots & \frac{\partial P_{\ell -1}}{\partial x_1} & \frac{\partial}{\partial x_1} \\
\vdots & \ddots & \vdots &\vdots  \\
\frac{\partial P_1}{\partial x_\ell} & \cdots & \frac{\partial P_{\ell -1}}{\partial x_{\ell }} & \frac{\partial}{\partial x_{\ell }} 
\end{array}
\right).
\end{displaymath}

Next we define an affine connection 
$\nabla :\mathrm{Der}_V\times \mathrm{Der}_V \rightarrow \mathrm{Der}_V$.

\begin{Def}\label{def:conn}
For given $\delta_1, \delta_2 \in \mathrm{Der}_V$ with 
$\delta_2 =\sum\limits_{i=1}^\ell f_i \frac{\partial}{\partial x_i}$, 
\begin{displaymath}
\nabla_{\delta_1}\delta_2 :=
\sum_{i=1}^\ell (\delta_1 f_i) \frac{\partial}{\partial x_i}.
\end{displaymath}
\hfill$\square$
\end{Def}

The connection $\nabla$ can be also characterized by the formula: 
\begin{equation}\label{eq:formula}
(\nabla_{\delta_1}\delta_2)\alpha =\delta_1(\delta_2\alpha),\ 
\forall \mbox{ linear function } \alpha \in V^\ast.
\end{equation}
This formula plays an important role in our computations. 

The derivation $\nabla_D$ by the primitive vector field is 
particularly important. Define 
$\mathbb{R}[V]^{W,\tau}:=\{f\in \mathbb{R}[V]^W |\ Df=0\}
=\mathbb{R}[P_1,\cdots ,P_{\ell -1}]$. Then $\nabla_D$ is an 
$\mathbb{R}[V]^{W,\tau}$-homomorphism. The following 
decomposition of 
$\mathrm{Der}_V^W={\rm D}^1(\mathcal{A})^W$ 
has been obtained in \cite{sai-linear},\cite{sai-lect}.

\begin{Thm}\label{thm:decomp}
Let $n\geq 1$, define 
\begin{displaymath}
\mathcal{G}_n:=\left\{ \delta \in \mathrm{Der}_V^W \left| \ (\nabla_D)^n\delta \in 
\sum_{i=1}^\ell \mathbb{R}[V]^{W,\tau}\frac{\partial}{\partial P_i}\right.\right\},
\end{displaymath}
then for every $n\geq 0$, 
$\nabla_D$ induces an $\mathbb{R}[V]^{W,\tau}$-isomorphism 
$\mathcal{G}_{n+1}\tilde{\rightarrow}\mathcal{G}_n$ 
and 
\begin{displaymath}
\mathrm{D}^1(\mathcal{A})^W=\bigoplus\limits_{n\geq1}\mathcal{G}_n.
\end{displaymath}
If we define $\mathcal{H}^k:=\bigoplus\limits_{n\geq k}\mathcal{G}_n$, 
then it becomes a rank $\ell$ free 
$\mathbb{R}[V]^W$-submodule of $\mathrm{Der}_V^W$, which is 
called the Hodge filtration.
\hfill$\square$
\end{Thm}

\noindent
In particular, 
$\nabla_D:\mathcal{H}^2\tilde{\rightarrow}\mathcal{H}^1=
\mathrm{D}^1(\mathcal{A})^W$ is an $\mathbb{R}[V]^{W,\tau}$-isomorphism. 
All we need in the sequel is the existence of an injection 
$\nabla_D^{-1}:\mathrm{Der}_V^W\rightarrow \mathrm{Der}_V^W$.

\section{Construction of a basis}\label{sec:two}

We construct a basis of $\mathrm{D} (\mathcal{A}^{(2k +\tilde{m})})$. 
The following is a key lemma which connects two filtrations, 
the Hodge filtration and the contact-order filtration.

\begin{Lemma}\label{lem:key}
Let $\delta' ,\delta \in \mathrm{Der}_V$ be vector fields on $V$ and 
assume $\nabla_D\delta'=\delta$. 
Then for any $H\in \mathcal{A}$, $\delta \alpha_H$ is divisible by 
$\alpha_H^m$ if and only if $\delta'\alpha_H$ is 
divisible by $\alpha_H^{m+2}$.
\end{Lemma}
{\bf Proof.}\\
Suppose $\delta '\alpha =\alpha^{m'}f$
(where $\alpha =\alpha_H$).
Then from (\ref{eq:formula}), 
\begin{equation}\label{eq:formula2}
(\nabla_D\delta ')\alpha =D(\delta '\alpha)=
\frac{1}{Q}
\det
\left(
\begin{array}{cccc}
\frac{\partial P_1}{\partial x_1} & \cdots & \frac{\partial P_{\ell -1}}{\partial x_1} & \frac{\partial}{\partial x_1} (\alpha^{m'}f)\\
\vdots & \ddots & \vdots &\vdots  \\
\frac{\partial P_1}{\partial x_\ell} & \cdots & \frac{\partial P_{\ell -1}}{\partial x_{\ell }} & \frac{\partial}{\partial x_{\ell }} (\alpha^{m'}f)
\end{array}
\right).
\end{equation}
Thus $\delta \alpha$ is divisible by $\alpha ^{m'-2}$. Further, 
assume $f$ is not divisible by $\alpha$, let us show that 
$\delta \alpha$ is not divisible by $\alpha ^{m'-1}$. Take a 
coordinate system $x_1,\cdots ,x_{\ell -1},x_\ell$ such that 
$x_\ell =\alpha$, then it suffices to show that 
\begin{displaymath}
\det \left(
\begin{array}{ccc}
\frac{\partial P_1}{\partial x_1}&\cdots &\frac{\partial P_{\ell -1}}{\partial x_1} \\
\vdots  &\ddots &\vdots \\
\frac{\partial P_1}{\partial x_{\ell-1}} &\cdots & \frac{\partial P_{\ell -1}}{\partial x_{\ell-1}}
\end{array}
\right)
\mbox{\ is not divisible by } \alpha.
\end{displaymath}
After taking $\otimes \mathbb{C}$ and restricting to 
$H_{\mathbb{C}}:=H\otimes \mathbb{C}$, 
determinant above can be interpreted as the Jacobian of the composed 
mapping 
\begin{displaymath}
\begin{array}{cccc}
\phi:&H_{\mathbb{C}} &\rightarrow & \mathrm{Spec}\mathbb{C}[V]^{W,\tau}\\
&&&\\
&(x_1,\cdots ,x_{\ell -1},0)&\mapsto&(P_1,\cdots ,P_{\ell -1}).
\end{array}
\end{displaymath}
On the other hand, since the set 
\begin{displaymath}
\{ x\in V\otimes \mathbb{C}|\ P_1(x)=\cdots =P_{\ell-1}(x)=0\}
\end{displaymath}
is a union of some eigenspaces of Coxeter transformations in $W$, 
which are regular, that is, they intersect with $H_{\mathbb{C}}$ 
only at $0\in H_{\mathbb{C}}$ \cite{sai-linear},\cite{sai-lect}.
Hence $\phi^{-1}(0)=\{0\}\subset H_{\mathbb{C}}$, 
and the Jacobian of $\phi$ cannot 
be identically zero. 
\hfill$\square$\\

\begin{Rem}\normalfont
The precise expression of the Jacobian of $\phi$ is 
obtained in \cite{sai-poly}. It is equal to 
the reduced defining equation of the union of hyperplanes 
$\bigcup\limits_{H'\in \mathcal{A}\backslash \{ H\}}(H\cap H')$, 
on $H$.
\end{Rem}

Because of Theorem \ref{thm:decomp} 
the operator $\nabla_D^{-1}$ is well defined on 
$\mathrm{Der}_V^W=\mathrm{D}^1(\mathcal{A})^W$, 
we have 

\begin{Lemma}\label{lem:prim}
Let $\delta \in \mathrm{Der}_V^W$ be a $W$-invariant vector field on $V$. 
Then for any $H\in \mathcal{A}$, $\delta \alpha_H$ is divisible by 
$\alpha_H^m$ if and only if $(\nabla_D^{-1}\delta)\alpha_H$ is 
divisible by $\alpha_H^{m+2}$.$\square$
\end{Lemma}

By induction with $\mathcal{H}^1=\mathrm{D}^1(\mathcal{A})^W$, 
Lemma \ref{lem:prim} indicates 
\begin{equation}\label{eq:filt}
\mathcal{H}^k=\nabla_D^{-k+1}\mathrm{D}^1(\mathcal{A})^W \subset 
\mathrm{D}^{2k+1}(\mathcal{A})^W.
\end{equation}
The converse is also true, which will be proved in 
\S\ref{sec:rem}.

We denote by 
$E:=\sum\limits_{i=1}^\ell x_i\frac{\partial}{\partial x_i}$
the Euler vector field. Note that $E$ is contained in 
$\mathrm{D}^1(\mathcal{A})^W$, $\nabla_\delta E=\delta$ and 
$\nabla_E \delta =(\deg \delta)\delta$ for any homogeneous 
vector field $\delta \in \mathrm{Der}_V$.
By Theorem \ref{thm:decomp}, we have a ``universal'' vector field 
$\nabla_D^{-k}E$.

As in \S\ref{sec:intro}, let $\tilde{m}:\mathcal{A}\rightarrow \{0,1\}$ 
be a multiplicity and assume that 
$\delta_1,\delta_2,\cdots,\delta_\ell \in \mathrm{D}(\mathcal{A}^{(\tilde{m})})$ 
be a free basis of the multiarrangement $\mathcal{A}^{(\tilde{m})}$. 

\begin{Thm}\label{thm:bas}
Under the above hypothesis, 
$\nabla_{\delta_1}\nabla_D^{-k}E, \cdots ,\nabla_{\delta_\ell}\nabla_D^{-k}E$
form a free basis of $\mathrm{D}(\mathcal{A}^{(\tilde{m}+2k)})$.
\end{Thm}
{\bf Proof.}\\
Let $\delta \in \mathrm{D}(\mathcal{A}^{(\tilde{m})})$, 
we first prove 
$\nabla_{\delta}\nabla_D^{-k}E \in \mathrm{D}(\mathcal{A}^{(\tilde{m}+2k)})$.
From (\ref{eq:filt}),  
$\nabla_D^{-k}E \in \mathrm{D}^{2k+1}(\mathcal{A})$, we may assume 
\begin{equation}\label{eq:fct}
(\nabla_D^{-k}E)\alpha =\alpha ^{2k+1} f
\end{equation}
for 
$\alpha =\alpha _H,\ (H\in \mathcal{A})$. 
Applying $\delta$ to the both sides of (\ref{eq:fct}), 
we have 
\begin{equation}
(\nabla_{\delta}\nabla_D^{-k}E)\alpha 
=\alpha ^{2k}\left( (2k+1)(\delta \alpha)f +\alpha (\delta f)\right).
\end{equation}
Since $\delta \alpha$ is divisible by $\alpha $ with multiplicity 
$\tilde{m}(H)\leq 1$, hence $(\nabla_{\delta}\nabla_D^{-k}E)\alpha$ is 
divisible by $\alpha ^{\tilde{m}(H) +2k}$.

Here we recall G. Ziegler's criterion on freeness of multiarrangements. 

\begin{Thm}\label{thm:zeg}{\normalfont \cite{zie}}
Let $\tilde{m}:\mathcal{A}\rightarrow \mathbb{Z}_{\geq 0}$ be a 
multiplicity and 
$\delta_1, \cdots ,\delta_\ell \in \mathrm{D}(\mathcal{A}^{(\tilde{m})})$ 
be homogeneous and linearly independent over $\mathbb{C}[V]$. 
Then $\mathcal{A}^{(\tilde{m})}$ is free with basis 
$\delta_1, \cdots ,\delta_\ell$ if and only if 
\begin{displaymath}
\sum_{i=1}^\ell \deg \delta_i =\sum_{H\in \mathcal{A}}\tilde{m}(H).
\end{displaymath}
\hfill$\square$
\end{Thm}

We compute the degrees of 
$\nabla_{\delta_1}\nabla_D^{-k}E, \cdots ,\nabla_{\delta_\ell}\nabla_D^{-k}E$, 
\begin{equation}\label{eq:degree}
\begin{array}{ccc}
\sum_{i=1}^{\ell} \deg (\nabla_{\delta_i}\nabla_D^{-k}E) 
&=& \sum_{i=1}^{\ell}\left( kh + \deg \delta _i \right) \\
&&\\
&=& kh\ell + \sum _{i=1}^{\ell}\deg \delta _i ,
\end{array}
\end{equation}
where $h=\deg P_\ell$ is the Coxeter number. On the other hand, 
the sum of multiplicities is 
\begin{equation}\label{eq:number}
\sum_{H\in \mathcal{A}}(\tilde{m}(H) +2k) =2k|\mathcal{A}|+
\sum_{H\in \mathcal{A}}\tilde{m}(H)
\end{equation}
The assumption implies 
$\sum_{H\in\mathcal{A}}\tilde{m}(H)=\sum_{i=1}^{\ell}\deg \delta_i$ and 
because of (\ref{eq:num}), 
we conclude that (\ref{eq:degree}) coincides with (\ref{eq:number}).
\hfill$\square$

\section{Some conclusions}\label{sec:rem}

\begin{Lemma}\label{lem:rel}  
$\ \nabla_{\frac{\partial}{\partial P_i}} \mathrm{D}^{2k+1}(\mathcal{A})^W
\subset \mathrm{D}^{2k-1}(\mathcal{A})^W\ (k>0)$.
\end{Lemma}
{\bf Proof.}\\
We only prove for $i=\ell$, remaining cases can be proved similarly. 
It is sufficient to show that $(\nabla_D\delta)\alpha_{H_0}$ has no poles 
for any 
$\delta \in \mathrm{D}^{2k+1}(\mathcal{A})^W$ and $H_0\in \mathcal{A}$. 
By (\ref{eq:formula2}), 
$QD\delta \alpha_{H_0}$ can be divided by $\alpha_{H_0}$, 
so all we have to show is that $QD\delta \alpha_{H_0}$ is divisible by 
$\beta :=\alpha _{H'}$ for all $H'\in \mathcal{A}\backslash \{H_0\}$.
We denote by $s_\beta \in W$ the reflection with respect to the hyperplane 
$H'\subset V$, then $s_\beta(\alpha)$ is expressed in the form 
$s_\beta(\alpha)=\alpha +2c\beta$ for some $c\in \mathbb{R}$. 
Apply $s_\beta$ to the function $QD\delta \alpha$, since $D$ and 
$\delta$ are  $W$-invariant, and $s_\beta (Q)=-Q$, 
$$
s_\beta (QD\delta \alpha ) = -QD\delta \alpha -2c QD\delta \beta .
$$
by using the equation $s_\beta(QD\delta \beta)=QD\delta \beta$, we have 
\begin{displaymath}
s_\beta (QD\delta \alpha + cQD\delta \beta)=
-(QD\delta \alpha + cQD\delta \beta).
\end{displaymath}
So $QD\delta \alpha + cQD\delta \beta$ is divisible by $\beta$, but from 
the first half of this proof, $cQD\delta \beta$ is divisible by $\beta$, 
and the other term $QD\delta \alpha$ is also divisible by $\beta$.
\hfill$\square$

As a consequence of induction, we have 
\begin{Cor}\normalfont{\cite{ter-hodge}} 
$\mathcal{H}^k =\mathrm{D}^{2k+1}(\mathcal{A})^W$.\hfill$\square$
\end{Cor}

Finally, we apply Theorem \ref{thm:main} to $\tilde{m}\equiv 0$ or 
$\tilde{m}\equiv 1$, since both $\mathrm{D}^0(\mathcal{A}) =\mathrm{Der}_V$ 
and $\mathrm{D}^1(\mathcal{A})$ are free, we obtain 

\begin{Cor}\normalfont{\cite{ter-multi}}
$\mathrm{D}^m(\mathcal{A})$ is free for all $m\geq 0$.\hfill$\square$
\end{Cor}

\noindent {\bf Acknowledgement.} I am grateful to Professor Kyoji Saito and Professor Hiroaki Terao for their useful comments.

\end{document}